\documentclass[11pt]{article}
\usepackage{amsmath}
\usepackage{bm}
\usepackage{psfrag}
\usepackage[all,cmtip,line]{xy}
\usepackage[normalem]{ulem}
\usepackage{tikz} 
\usepackage[utf8]{inputenc}
\usepackage{pgfplots}
\pgfplotsset{compat=1.7}
\usepackage{scalerel}
\usepackage{amssymb}
\usepackage{caption}
\usepackage{subcaption}
\usepackage[shortlabels]{enumitem}
\usepackage{relsize}
\usepackage{authblk}
\usepackage[margin=01.0in]{geometry}
\usepackage{hyphenat}

\newcommand{\IZ}{\mathbb{Z}}
\newcommand{\IR}{\mathbb{R}}
\newcommand{\IN}{\mathbb{N}}

\newcommand{\cA}{\mathcal{A}}
\newcommand{\cB}{\mathcal{B}}
\newcommand{\cN}{\mathcal{N}}

\newcommand{\cV}{\mathcal{V}}

\newcommand{\cO}{\mathcal{O}}
\newcommand{\cF}{\mathcal{F}}

\newcommand{\p}[1]{\left\langle #1 \right\rangle}
\newcommand{\norm}[1]{\Vert #1 \Vert}

\begin{document}

\title{Optimization methods for achieving high diffraction efficiency with perfect electric conducting gratings}

\author[1]{Rub\'en Aylwin}
\author[2]{Gerardo Silva-Oelker}
\author[3]{Carlos Jerez-Hanckes}
\author[4]{Patrick Fay}
{\tiny 
\affil[1]{Department of Electrical Engineering, Pontificia Universidad Cat\'olica de Chile, Santiago, Chile.}
\affil[2]{Department of Mechanical Engineering, Universidad Tecnol\'ogica Metropolitana, Santiago, Chile.}
\affil[3]{Faculty of Engineering and Sciences, Universidad Adolfo Ib\'a\~nez, Santiago, Chile.}
\affil[4]{Department of Electrical Engineering, University of Notre Dame, IN, USA.}
}
\maketitle




\begin{abstract}
This work presents the implementation, numerical examples and
experimental convergence study of first- and second-order optimization methods applied to one-dimensional periodic gratings. Through boundary integral equations and shape derivatives, the profile of a grating is optimized such that it maximizes the diffraction efficiency for given diffraction modes for transverse electric polarization. We provide a thorough comparison of three different optimization methods: a first-order method (gradient descent); a second-order approach based on a Newton iteration, where the usual Newton step is replaced by taking the absolute value of the eigenvalues given by the spectral decomposition of the Hessian matrix to deal with non-convexity;
and the Broyden-Fletcher-Goldfarb-Shanno (BFGS) algorithm, a quasi-Newton
method. Numerical examples are provided to validate our claims. Moreover, two
grating profiles are designed for high efficiency in the Littrow configuration and then compared to a high efficiency commercial grating. Conclusions and recommendations, derived from the numerical experiments, are provided as well as future research avenues.
\end{abstract}

%
%
%
%
\section{Introduction}
\label{sec:intro}

One-dimensional gratings are able to diffract, split, reflect and transmit light, depending on geometrical parameters such as amplitude, period, and shape. For example, it is well known that metallic gratings exhibit Wood's anomalies \cite{Wood1902, Rayleigh1907} and cavity resonances \cite{Quaranta2018}. Of particular interest in applications is the diffraction efficiency---defined as the amount of power diffracted in one mode \cite{Petit1980}. Due to these properties, gratings have several applications in science and engineering, ranging from X-ray spectroscopy \cite{McEntaffer2013}, energy conversion devices such as photovoltaics \cite{Solano2013} and thermophotovoltaics \cite{SilvaOelker2018b}, beam splitters \cite{Davis2001}, quantum cascade lasers \cite{Dirisu2007}, and filters \cite{Niraula2014}, to name a few. Consequently, a large body of literature has been devoted to the study of physical phenomena \cite{Liang2013,Graham2014,Maradudin2016}, mathematical modeling \cite{Chandler-Wilde1999, Zhang2003, Hu2015,SilvaOelker2018}, fabrication \cite{Lu2009,Grigorescu2009, Saleem2014}, and optimization of gratings. 

Among the optimization methods for designing gratings, genetic algorithms and particle swarm optimization have been studied and implemented for applications such as energy conversion \cite{SilvaOelker2019, Chen2010} and filters \cite{ShokoohSaremi2007}. These methods have shown to be practical options due to their simplicity and flexibility in implementation. However, the aforementioned techniques suffer from disadvantages such as partial convergence and the need for a large number of evaluations of the cost function. In contrast, methods that rely on following the gradient direction to minimize a suitable cost function can be more appropriate for some problems. Roger \cite{Roger1983}, for instance, studied the optimization of a perfect electric conducting (PEC) grating by the steepest-descent and conjugate-gradient methods along with boundary integral equations (BIEs). Recently, Bao et al. \cite{Bao2013} implemented the Landweber iteration also along with BIE for inverse problems, showing the applicability of gradient optimization. Methods based on so-called {\em shape derivatives} have been theoretically proposed \cite{Eppler2006, Eppler2012} for optimization and successfully applied to the design of devices \cite{Paganini2015}. However, grating geometry optimization using shape derivatives, to the best of our knowledge, has not been explored yet. This work aims to provide details and examples of first- and second-order methods based on shape gradients for PEC grating profile optimization.

The use of optimization algorithms can lead to non-trivial grating shapes (e.g., \'echelle or holographic). Fortunately, mainstream fabrication approaches for manufacturing complex profiles have become available, including electron beam lithography \cite{Zeitner2012,Sudheer2016}, laser ablation \cite{Bader2006}, deposition-and-etch based approaches \cite{VilaComamala2018,Aryal2012}, and interference lithography \cite{khan:2016, gao:2019}.  While not as well developed commercially, intricate structures can also be fabricated using techniques such as nanoimprint  \cite{Yu2003,Honma2016}, surface wrinkling \cite{Wei2014}, and scanning-probe lithography \cite{lassaline2019optical}. With these tools, complex grating structures can now be manufactured, thereby increasing the range of shape possibilities and making it important to design and fabricate optimal gratings.

In this work, we implement and compare first- and second-order optimization algorithms to maximize the diffraction efficiency for transverse electric polarization. Diffraction efficiency was chosen as the target figure of merit for its importance in applications \cite{Bunkowski2006}. In our approach, the wave scattering model problem in computational volume is reduced to one defined on the grating boundary through an integral formulation and a suitable quasi-periodic Green's function. Optimization is carried out by minimizing or maximizing a cost or objective function using first- and second-order shape gradients \cite{Costabel2012} of the defined far-field operator, which maps the grating profile to the diffraction field components. In the case of the first-order method, it is shown that the shape derivative approach is equivalent to the work presented by Roger \cite{Roger1983}. However, since this optimization problem can be highly non-convex, additional techniques are explored to improve performance.  In this regard, two Newton-based second-order algorithms are implemented and studied. First, a second-order method with a modified step \cite{dauphin2014identifying, Paternain2019},able to deal with non-convexity, is considered.
Secondly, the quasi-Newton Broyden-Fletcher-Goldfarb-Shanno (BFGS)  algorithm is implemented for comparison. Numerical results show that, for certain objective functions, the modified Newton method converges quadratically, whereas the BFGS algorithm converges super-linearly. Moreover, these methods reduce both computational time and the number of iterations required to find a maximum. These results pave the way for more efficient approaches of grating optimization and periodic structures in general. Furthermore, the studied techniques can also be applied to inverse design problems.

The rest of the paper is organized as follows. In Section \ref{sec:bem}, the wave scattering by a grating problem is presented, along with the boundary integral formulation and the definition of diffraction efficiency. Section \ref{sec:opti} describes the optimization methods: shape derivatives are introduced and first- and second-order methods are detailed. Numerical examples are provided in Section \ref{sec:numexam}, analyzing computation times as well as designing grating profiles. Finally, conclusions are drawn in Section {\ref{sec:conc}}.

%
%
\section{Scattering Problem and Boundary Integral Formulation}
\label{sec:bem}
\begin{figure}[ht!]
	\centering\includegraphics[width=8cm]{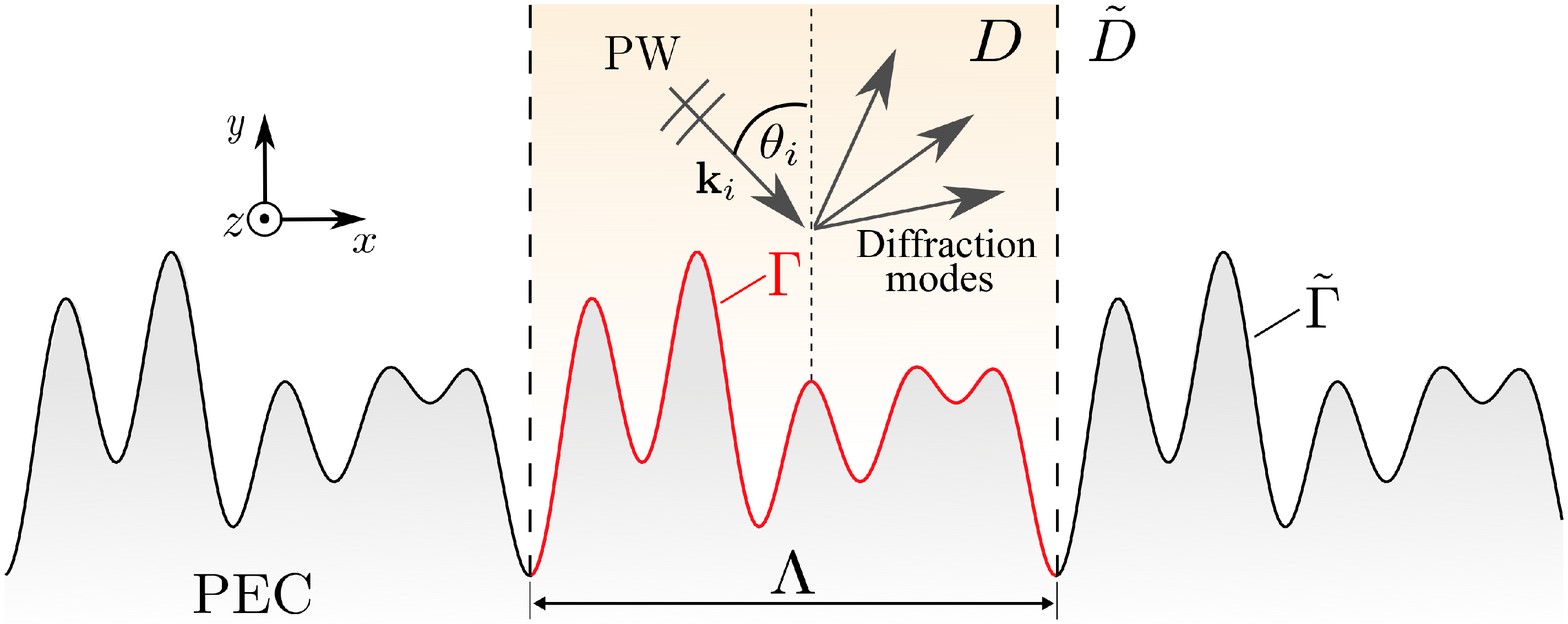}
	\caption{Perfectly conducting grating with a period $\Lambda$: the domain $\widetilde{D}$ above the grating and its surface $\widetilde{\Gamma}$ are shown. Both $D$ and $\Gamma$ are defined in one period.}
	\label{fig:grating}
\end{figure}

A PEC grating with period $\Lambda$ and surface given by $\widetilde{\Gamma}\times\IR$ is considered (see Fig.~\ref{fig:grating}), where $\widetilde{\Gamma}$ is assumed to be (at least) Lipschitz continuous. The domain $\widetilde{D}\times \mathbb{R}\subset \mathbb{R}^3$ is defined as the open region of propagation above the grating surface, assumed to be free space and characterized by its impedance $\eta = \sqrt{\mu/\varepsilon}$, where $\varepsilon$ and $\mu$ correspond to vacuum permittivity and permeability, respectively.

\subsection{Scattering Problem}
\label{ssec:scatterProb}
We consider the scattering of a monochromatic plane wave (PW) with wave vector $\mathbf{k}_i$ by a grating surface described in Fig.~\ref{fig:grating}. Furthermore, we assume the wave vector to be such that it has no $z$ component ($\mathbf{k}_i=(k_x,k_y)$), so that Maxwell's equations may be decomposed into transverse electric (TE) and transverse magnetic (TM) polarizations (see \cite{bonnet1994guided} and references therein). For brevity, only TE polarization ($E_{x}=E_{y}=H_{z}=0$) is considered in this study. Nevertheless, optimization for TM polarization follows the same approach with different boundary conditions. Hence, this work focuses on the $z$ component of the electric field, defined as:
\begin{equation}
	u(\mathbf{r}) = E_{z}(\mathbf{r}) \quad \mathbf{r} \in \widetilde{D},
	\label{eq:e_field_definition}
\end{equation} 
where $\mathbf{r}:=(x,y)\in\IR^2$ is the position vector. Then, the total field $u$ satisfies the Helmholtz equation:
\begin{align}
	\Delta u + k^{2}u &= 0 \quad \text{in } \quad \widetilde{D},
	\label{eq:helmholtz_whole_domain}\\
	u &= 0 \quad \text{on} \quad \widetilde{\Gamma},
\end{align}
where $k:=\omega/c$ is the wavenumber. As previously stated, we consider the scattering of an incident PW (see Fig.~\ref{fig:grating}): $$u^{i}(\mathbf{r}) := E_{0}e^{\imath \mathbf{k}_{i}\cdot \mathbf{r} }\quad \forall\; \mathbf{r}\in\widetilde{D},$$ onto the grating surface $\widetilde\Gamma$, where $\mathbf{k}_{i}=(k_{x},k_{y}) = k(\sin \theta_{i}, -\cos \theta_{i})$ is the incident wave vector with incidence angle $\theta_{i}$ and null $z$ component, and $E_0>0$ corresponds to the incident electric field amplitude. By linearity, the total field is $u = u^{i} + u^{s}$, where $u^{s}$ is the scattered field. Since PWs satisfy \eqref{eq:helmholtz_whole_domain}, we may restate the previous system as a problem for the scattered field as follows:
\begin{eqnarray}
	\Delta u^{s} + k^{2}u^{s} & = & 0 \qquad \text{in} \quad \widetilde{D},  \label{eq:helmholtz}\\
	u^{s} & = & -u^{i} \quad \text{on} \quad \widetilde{\Gamma}.
\end{eqnarray}
The grating periodicity and the quasi-periodicity of $u^i$, i.e.
\begin{align}
	u^i(x+n\Lambda,y)=u^i(x,y)e^{\imath k_x\Lambda n}\quad \forall\; n\in\IZ, \label{eq:quasiperiodicity}
\end{align}
enable us to restate the problem for the scattered field on one cell $D$---bounded in the $x$-direction---of the infinite periodic domain $\widetilde D$ by enforcing the condition in \eqref{eq:quasiperiodicity} at the left and right boundaries of $D$. However, the problem still lacks appropriate radiation conditions at infinity, which, for the sake of brevity, will not be  derived here. Instead, we refer to \cite{bonnet1994guided, nedelec1991integral, kirsch1993diffraction} and references therein for a derivation and discussion of adequate radiation conditions in periodic problems. For a fixed $d\in\IR$, such that $\widetilde\Gamma$ lies strictly below the line $y=d$,
we define the following domains:
	\begin{eqnarray*}
		D := \{\mathbf{r}\in\widetilde{D}\; :\; x\in(0,\Lambda)\} \quad \text{and}\quad 
		D_d :=\{\mathbf{r}\in{D}\; :\; y<d\}.
	\end{eqnarray*}
	 Analogously, we define the grating surface over one period
	\begin{equation*}
		\Gamma :=\{\mathbf{r}\in\widetilde{\Gamma}\; :\; x\in(0,\Lambda)\}.
	\end{equation*}
Additionally, the Sobolev space of square integrable functions with square integrable first-order derivatives on $D_d$ that satisfy \eqref{eq:quasiperiodicity} is denoted $H^1_{k_x;\Lambda}(D_d)$. We then seek to find a field $u^{s}\in H^1_{k_x;\Lambda}(D_d)$ that solves
\begin{eqnarray}
\begin{aligned}
	\Delta u^{s} + k^{2}u^{s} & =  0 \quad \text{in} \; D_d, \\
	u^{s} & =  -u^{i}=:g \quad \text{on }  \Gamma,\\
	u^{s}(x,y) & =  \sum_{n \in \mathbb{Z}}u_{n}(d) e^{\imath (k_{y,n}(y-d)+k_{x,n}x)} \quad \forall\; y\geq d, 
	\end{aligned}\label{eq:periodic_full_problem}
\end{eqnarray}
where,  for each $n\in\IZ$, $u_{n}(d)$ are the Rayleigh coefficients,  $k_{x,n}:=k \sin \theta_{i}+2\pi n  / \Lambda$, and
\begin{gather*}
	k_{y,n}:=\begin{cases}
		\sqrt{k^{2}-k_{x, n}^{2}}\quad &\text{if}\quad k^{2} \geq k_{x, n}^{2}\\
		\imath \sqrt{k_{x, n}^{2}-k^{2}}\quad &\text{if}\quad k^{2}<k_{x, n}^{2}
	\end{cases}.
\end{gather*}

The last line in \eqref{eq:periodic_full_problem} corresponds to the radiation condition for the periodic domain, also known as the Rayleigh expansion for $u^s$ \cite{bonnet1994guided, Petit1980, kirsch1993diffraction}. In order to avoid Rayleigh anomalies, we further assume that no $n\in\IZ$ exists such that $k_{x,n}^2=k^2$.  

\subsection{Boundary Integral Equations}
\label{ssec:boundary_integral_eqs}

To find the scattered field, we use BIEs---particularly the electric field integral equation---motivated by the unboundedness of $D$. Therefore, we make the following representation Ansatz:
\begin{equation}
	u^{s}(\mathbf{r}) = (\mathsf{SL} j)(\mathbf{r}), \quad \mathbf{r} \in D,
	\label{eq:indirect_formulation}
\end{equation} 
where $j$ is an unknown surface density (in this case, the surface current density) and $\mathsf{SL}$ is the quasi-periodic single-layer operator, whose action on the density $j$ may be represented through the quasi periodic Green's function $G^{p}$ as:
\begin{equation}
	(\mathsf{SL} j)(\mathbf{r}) := \imath k\eta\int_{\Gamma}G^{p}(\mathbf{r}, \mathbf{r}')j(\mathbf{r}')d\Gamma(\mathbf{r}'), \quad \mathbf{r} \in D.
	\label{eq:single_layer_potential}
\end{equation}

Using \eqref{eq:indirect_formulation}, \eqref{eq:single_layer_potential}, and the boundary conditions in  \eqref{eq:periodic_full_problem}, the BIE for the unknown surface density $j$ may be written as:
\begin{equation}
	\cV j = -u^{i} \quad \text{on } \Gamma,\label{eq:BIE}
\end{equation}
where $\cV$ is the boundary integral operator mapping the surface density $j$ to the boundary values of $\mathsf{SL} j$.
The quasi-periodic Green's function, $G^{p}$, is given by \cite{Linton1998, Tsang2001}:
\begin{equation}
	G^{p}(\mathbf{r},\mathbf{r}')  = \frac{\imath}{2\Lambda} \sum_{n=-\infty}^{\infty} \frac{e^{\imath k_{x,n} (x-x') + \imath k_{y,n}|y-y'|}}{k_{y,n}}.
	\label{eq:quasi-periodic_green_function}
\end{equation}

%
%
\subsection{The Scattered Far-Field}
\label{ssec:sffield}
As displayed in \eqref{eq:periodic_full_problem}, the scattered field may be decomposed as a linear combination of complex exponentials. Of these, only those corresponding to $n\in\IZ$ such that $k^2> k_{x,n}^2$ do not decay as $y$ grows to infinity---as noted above, we assume no $n\in\IZ$ exists such that $k^2 = k_{x,n}^2$. Therefore, the scattered far-field---the scattered field for $d$ far
from the grating surface---depends solely on the Rayleigh coefficients corresponding to a real value for $k_{y,n}$, i.e., on $\{u_n(d)\}_{n\in\cN({k,\theta_i,\Lambda})}$ where:
\begin{align*}
	\cN({k,\theta_i,\Lambda}):=\{n\in\IZ\ :\ k > \vert k\sin(\theta_i)+2\pi n/\Lambda\vert\}.
\end{align*}
The coefficients of the Rayleigh expansion in the radiation condition of Eq. (\ref{eq:periodic_full_problem}) may be extracted from the scattered field as follows:
\begin{align}
\begin{aligned}
	u_{n}(d)&=\frac{1}{\Lambda}\int_{0}^{\Lambda}u^s(x,d)e^{-\imath k_{x,n}x}dx\\&=-e^{\imath k_{y,n}d}\frac{k\eta}{2\Lambda k_{y,n}}\int_{\Gamma}j(\mathbf{r}){e^{-\imath\mathbf{K}_n\cdot\mathbf{r}}}d\Gamma(\mathbf{r}),
	\end{aligned}
	\label{eq:BEMffield}
\end{align}
where $\mathbf{K}_{n}:=(k_{x.n},k_{y,n})$. Then, the associated grating diffraction efficiency can be computed through $e_n:=\frac{k_{y,n}}{k_y}\vert u_n(d)\vert^2$, further justifying the choice to employ BIE to solve the scattering problem. Henceforth, we omit the dependency
of the coefficients $u_n$ on $d$, since it is of no significance for the computation
of the efficiencies, as indicated by \eqref{eq:BEMffield} (the magnitude $-u_n(d)e^{-\imath k_{y,n}d}$ is constant in $d$). Coefficients $u_n$ are thus calculated using
\begin{align}
u_{n}:=\frac{k\eta}{2\Lambda k_{y,n}}\int\limits_{\Gamma}j(\mathbf{r}){e^{-\imath\mathbf{K}_n\cdot\mathbf{r}}}d\Gamma(\mathbf{r}).\label{eq:BEMffield_ind_d}
\end{align}

%
%
\section{Optimization}
\label{sec:opti}
We aim to find optimal grating profiles
in the sense that they maximize (or minimize) functions of the diffraction efficiencies, e.g.,
\begin{align*}
	\max\limits_{\Gamma}e_n(\Gamma),\quad \min\limits_{\Gamma}\sum\limits_{n\in \cN'}(e_n(\Gamma)-e^{\text{obj}}_{n})^2,
\end{align*}
for some $\cN'\subseteq\cN(k,\theta_i,\Lambda)$ and some objective efficiencies $e^{\text{obj}}_n$  in the case of minimization. With this in mind, we turn to the problem of computing derivatives of the far-field
Rayleigh coefficients as functions of the grating geometry  for the implementation of local search optimization algorithms. To do so, we introduce the $n$-th mode far-field operator as
\begin{align}\label{eq:ffieldop}
	\cF_n:
	\Gamma \mapsto u_n.
\end{align}
This operator maps a grating profile $\Gamma$ to $\cF_n(\Gamma)=u_n$, the $n$-th coefficient
of the Rayleigh expansion of the diffracted field.
Henceforth, we assume $k$, $\Lambda$, and $\theta_i$ to be fixed.
\subsection{Shape Calculus and Shape Derivatives}\label{ssec:shapeder}
To compute derivatives of the far-field operator with respect to the grating geometry, we use tools from shape calculus \cite{Pironneau1984,sokolowski1992introduction}. For $m\in\IN$ and a given grating profile $\Gamma$ of class $C^m$ (the $m$-th derivative exists and is continuous) we introduce the following space of periodic functions: $$C^\ell_{p}(\Gamma):=\left\{\tau\in C^\ell(\Gamma;\IR^2)\ :\ \tau \ \text{is periodic }\right\},$$ 
for any non-negative integer $\ell\leq m$. Moreover, for any $a\in C_p^1(\Gamma)$, we introduce the perturbed grating profile as
$$\Gamma_a:=\{\mathbf{r}+a(\mathbf{r})\ :\ \mathbf{r}\in\Gamma \},$$
whenever $a$ is small enough so that $\Gamma_a$ does not self intersect. Then, the definition of the shape derivative of $\cF_n$ at $\Gamma$ in the direction $a$ is
\begin{align}
	\cF'_n(\Gamma;a):=\lim\limits_{\epsilon\rightarrow 0}\frac{1}{\epsilon}(\cF_n(\Gamma_{\epsilon a})-\cF_n(\Gamma)),\label{eq:firsDer}
\end{align}
whenever the limit in the right-hand side of \eqref{eq:firsDer} exists \cite{kirsch1993domain, kirsch1993diffraction, Pironneau1984, sokolowski1992introduction}. Moreover, we say $\cF_n$ is shape-differentiable on $C^\ell_p(\Gamma)$ if $\cF_n'(\Gamma;\cdot)$ is a linear and continuous functional on $C^\ell_p(\Gamma)$.
The computation of $\cF'_n(\Gamma;a)$,
for $a$ in $C^2_p(\Gamma)$, can be performed by solving a problem analogous to \eqref{eq:periodic_full_problem} through the use of the following relation:
\begin{align*}
	\cF'_n(\Gamma;a)=\cF_n(u'[a]),
\end{align*}
where $u'[a]$ is the unique solution in $H^1_{k_x;\Lambda}(D_d)$ to the following scattering problem:
	\begin{align}
	\begin{gathered}
		\Delta u'[a] + k^{2}u'[a]  = 0 \quad \text{in} \; D_d, \\
		u'[a]  =  -(a\cdot\nu)\frac{\partial}{\partial \nu}u \quad \text{on }  \Gamma,\\
		u'[a](\mathbf{r})  = \sum_{n \in \mathbb{Z}}u'_{n}[a](d){e^{\imath (k_{y,n}(y-d)+k_{x,n}x)}} \quad \forall\; y\geq {d},
		\end{gathered}\label{prob:first_der}
	\end{align}
	where $\nu$ is the normal vector to $\Gamma$. This result follows by
	modifying the proof of Theorem 2.1 in \cite{kirsch1993domain} to
	the quasi-periodic setting (\emph{cf.}~\cite{APJ:2019}).

The second-order shape derivative of $\cF_n$ at $\Gamma$ along the directions $a_1$, $a_2 \in C^1_{p}(\Gamma)$ is defined as
\begin{align}
	\cF_n''(\Gamma;a_1;a_2):=\lim\limits_{\epsilon\rightarrow 0}\frac{1}{\epsilon}(\cF_n'(\Gamma_{\epsilon a_2};a_1\circ\psi_\epsilon)-\cF_n'(\Gamma;a_1)),\label{eq:secder}
\end{align}
with 
\begin{align*}
	\psi_\epsilon:=\phi_{\epsilon}^{-1},\quad \phi_{\epsilon}:=\begin{cases}\Gamma &\to\Gamma_{\epsilon a_2}\\
		\mathbf{r} &\mapsto \mathbf{r}+\epsilon a_2(\mathbf{r})
	\end{cases}.
\end{align*}
This definition was introduced in the context of scattering from bounded obstacles in \cite{hettlich1999second} and satisfies $\cF_n''(\Gamma;a_1;a_2)=\cF_n''(\Gamma;a_2;a_1)$, i.e., it is symmetric. Analogous to the case of first-order shape derivatives, the computation of $\cF_n''(\Gamma;a_1;a_2)$
for $a_1$ and $a_2$ in $C^3_p(\Gamma)$ follows from:
	\begin{align*}
		\cF_n''(\Gamma;a_1;a_2)=\cF_n(u''[a_1;a_2]),
	\end{align*}
	where $u''[a_1;a_2]$ is the unique solution in $H^1_{k_x;\Lambda}(D_d)$ to the following
	scattering problem:
	\begin{gather*}
		\Delta u''[a_1,a_2] + k^{2}u''[a_1,a_2]  = 0 \quad \text{in} \; D_d, \\
		 u''[a_1 ,a_2] = -(a_1\cdot\nu)\frac{\partial u'[a_2]}{\partial \nu} - (a_2\cdot\nu)\frac{\partial u'[a_1]}{\partial \nu}\\
		+ \mathlarger{\mathlarger{(}}(a_1\cdot\nu)(a_2\cdot\nu)-(a_1\cdot\tau)(a_2\cdot\tau)\mathlarger{\mathlarger{)}}\kappa\frac{\partial u}{\partial \nu}\\
		 + \mathlarger{\mathlarger{(}} (a_1\cdot\tau)(\tau\cdot\nabla(a_2\cdot\nu))+(a_2\cdot\tau)(\tau\cdot\nabla(a_1\cdot\nu)) \mathlarger{\mathlarger{)}}\frac{\partial u}{\partial \nu}
		\quad \text{on }  \Gamma,\\
		 u''[a_1,a_2](\mathbf{r}) = \sum_{n \in \mathbb{Z}}u''_{n}[a_1;a_2](d){e^{\imath (k_{y,n}(y-d)+k_{x,n}x)}} \quad \forall\; y\geq {d},
	\end{gather*}
	where $\kappa$ denotes the curvature of $\Gamma$, and $\nu$ and $\tau$ denote the normal and tangent vectors to the surface, respectively.

Throughout the following sections, we consider, without loss of generality, $\Lambda\equiv1$. Moreover, recalling $\mathbf{r}=(x,y)\in\IR^2$
 we assume $\Gamma$ to be given as follows 
\begin{align*}
	\Gamma:=\left\lbrace\mathbf{r}\ :\ y=\sum\limits_{\ell=1}^{N}\cA_\ell\sin(2\pi \ell x)+\cB_\ell\cos(2\pi \ell x),\ x\in (0,1)\right\rbrace,
\end{align*}
thus fulfilling the conditions that allow us to compute first and second-order shape derivatives as solutions to the aforementioned boundary value
problems. Then, $\cF_n(\Gamma)=\cF_n(\cA_1,\cB_1,\hdots ,\cA_N, \cB_N)$ and the optimization problem is now finite-dimensional over the $2N$ variables $\{\cA_\ell\}_{\ell=1}^{N}$ and $\{\cB_\ell\}_{\ell=1}^{N}$. It holds that
\begin{align*}
	\frac{\partial\cF_n}{\partial\cA_\ell}
	=\cF_n'(\Gamma;{\mathsf{e}_y} \sin(2\pi \ell x))\quad \text{and} \quad
	\frac{\partial\cF_n}{\partial\cB_\ell}
	=\cF_n'(\Gamma;{\mathsf{e}_y }\cos(2\pi \ell x)),
\end{align*}
where $\mathsf{e}_y$ is the canonical vector in the $y$-direction.
Note that $\mathsf{e}_y \sin(2\pi\ell x)$ and $\mathsf{e}_y \cos(2\pi\ell x)$ represent directions, i.e., we need only consider perturbations
of the form $a(\bm{r})=\mathsf{e}_y a_y(\bm{r})$, with $a_y$ scalar.
Moreover, our choice of $\Gamma$ implies that $2N$ problems are solved to compute every required derivative. Analogous relations can be found for the second-order derivatives of $\cF_n$. 

%
%
\subsection{Adjoint Method for the Computation of Shape Derivatives}
The computation of each derivative of
$\cF_n(\cA_1,\cB_1,\hdots ,\cA_N, \cB_N)$ requires solving a
BIE with a varying right-hand side.
The same holds for each second-order derivative, which also requires the computation of boundary data of first-order shape derivatives. We may achieve the computation of derivatives more efficiently
through the adjoint method \cite{Gstrang:2007}. 

As noted previously, $j$ is the solution to \eqref{eq:BIE}, which also satisfies $j=\frac{\partial}{\partial\nu}(u^s+u^i)$. Let us also introduce $\tilde{j}$ as a surface density such that $u'[a](\mathbf{r})=(\mathsf{SL}\tilde j)(\mathbf{r})$ for $\mathbf{r}$ in $D$. Since the shape derivative, $u'[a]$, satisfies \eqref{prob:first_der}, then $\tilde j$ solves a BIE similar to \eqref{eq:BIE}, but considering the boundary condition in \eqref{prob:first_der}, i.e., $\cV \tilde j = -(a\cdot\nu)j$ on $\Gamma$. Recognizing that the coefficients $u_{n}$ in \eqref{eq:BEMffield_ind_d} can be written as a duality product and denoting $g_n(\mathbf{r}):=e^{-i\mathbf{K}_n\cdot\mathbf{r}}$, we obtain for $a\in C^2_p(\Gamma)$:
\begin{eqnarray*}
	\cF_n'(\Gamma;a) &=& \frac{k\eta}{2\Lambda k_{y,n}}\p{\tilde j,g_n}_{\Gamma}=\frac{k\eta}{2\Lambda k_{y,n}}\p{\cV^{-1}(-(a\cdot\nu)j),g_n}_{\Gamma}\\
	&=& \frac{k\eta}{2\Lambda k_{y,n}}\p{\cV^{-\top}(g_n),-(a\cdot\nu)j}_{\Gamma}.
\end{eqnarray*}

Hence, if $j^{\text{adj}}$ is a surface density that solves the BIE:
\begin{align}
	\cV^{\top}j^{\text{adj}}=g_n,\label{eq:V_adj_bie}
\end{align}
then, we may compute shape derivatives in the direction $a\in C^2_p(\Gamma)$ as:
\begin{align}
	\cF_n'(\Gamma;a)=\frac{k\eta}{2\Lambda k_{y,n}}\p{j^{\text{adj}},-(a\cdot\nu)j}_{\Gamma}.\label{eq:shorthandderivative}
\end{align}
According to (\ref{eq:shorthandderivative}), all of the first-order derivatives can be computed by solving two BIEs of the form of \eqref{eq:BIE} and then computing the specified integral. In \eqref{eq:V_adj_bie}, $\cV^\top$ is the boundary integral operator analogous to $\cV$ for the adjoint quasi-periodic Green's function given by
\begin{align*}
	G^p_{\text{adj}}(\mathbf{r},\mathbf{r}') = \frac{\imath}{2\Lambda} \sum_{n=-\infty}^{\infty} \frac{e^{\imath \widetilde{k}_{x,n} (x-x') + \imath k_{y,n}|y-y'|}}{k_{y,n}},
\end{align*}
where $\widetilde{k}_{x,n}:=-k\sin(\theta_i)+n$ for all $n\in\IN$.
The same
formula for the first-order shape derivatives of the far-field operator was found in \cite{roger1982generalized}
through the use of reciprocity relations between solutions of scattering problems.
Analogously, 
\begin{align}
	\cF_n''(\Gamma;a_1;a_2)=\frac{k\eta}{2\Lambda k_{y,n}}\p{j^{\text{adj}},u''[a_1,a_2]\vert_{\Gamma}}_{\Gamma}.\label{eq:shorthandsecderivative}
\end{align}

%
%
\subsection{First-Order Optimization}
\label{ssec:firstopti}
Throughout this and the following two subsections, we consider the optimization problem of minimizing a function $f:\IR^{2N}\to\IR$ such that 
\begin{align*}
	f(\cA_1,\cB_1,\hdots,\cA_N,\cB_N)=\tilde{f}(e_n(\cA_1,\cB_1,\hdots,\cA_N,\cB_N))
\end{align*}
for a smooth function $\tilde{f}:\IR\to\IR$ (e.g., $\tilde{f}(x)=x$, $\tilde{f}(x)=-x$ or $\tilde{f}(x)=(x-c)^2$ for some $c\in\IR$). Hence, our optimization problem may be stated as
\begin{align*}
	\min\limits_{\mathbf{x}\in\cO_{\text{ad}}} f(\mathbf{x})=\min\limits_{\Gamma\in\cO^{S}_{\text{ad}}}\tilde{f}(e_n(\Gamma)),
\end{align*}
where $\cO_{\text{ad}}$ is a set of admissible shape parameters and $\cO^S_{\text{ad}}$ is the set of admissible boundaries $\Gamma$ such that each $\mathbf{x}\in\cO_{\text{ad}}$ determines exactly one $\Gamma\in\cO^S_{\text{ad}}$. Notice that first-order shape derivatives of the grating efficiency may be computed as
\begin{align}
	\begin{aligned}
		e_n'(\Gamma;a)=\frac{2k_{y,n}}{k_y}\left(\text{Re}(\cF_n(\Gamma))\text{Re}(\cF_n'(\Gamma;a))\right.\\
		+\left.\text{Im}(\cF_n(\Gamma))\text{Im}(\cF_n'(\Gamma;a))\right).
	\end{aligned}\label{eq:efficiencyder}
\end{align}
Given \eqref{eq:shorthandderivative} and \eqref{eq:efficiencyder}, first-order derivatives of
the objective function can be computed, requiring the solution of only two integral equations
and $2N$ independent integrals, which may be computed in parallel. 

For an initial set-up of the
geometry given by $\IR^{2N}\ni \bm{x}^{(0)}:=(\cA_1^{(0)},\cB_1^{(0)},\hdots,\cA_N^{(0)},\cB_N^{(0)})$,
we will consider the usual first-order optimization algorithm (steepest descent):
\begin{align}
	\bm{x}^{(t)}=\bm{x}^{(t-1)}-h^{(t-1)}\nabla f(\bm{x}^{(t-1)}),
	\label{eq:graddescent}
\end{align}
where $h^{(t-1)}>0$ is the step size of the method at iteration $t$. The step size at each iteration may be found exactly, as
$$h^{(t-1)}=\arg\min\limits_{h>0}f(\bm{x}^{(t-1)}-h\nabla f(\bm{x}^{(t-1)})),$$
or by employing a backtracking strategy. Since choosing the
step size exactly is impractical, we choose to backtrack the step size through the Armijo-Goldstein rule: for a fixed $\alpha\in (0,1)$, $\beta\in (0,\frac{1}{2})$ and initial step size estimate $h>0$, we check whether the Armijo-Goldstein condition is satisfied:
\begin{align}
f(\bm{x}^{(t-1)}+h\bm{p}^{t-1})<f(\bm{x}^{(t-1)})-\beta h\bm{p}^{t-1}\cdot\nabla f(\bm{x}^{(t-1)}),\label{eq:armijo_rule}
\end{align}
where $\bm{p}^{(t-1)}=\nabla f(\bm{x}^{(t-1)})$. If \eqref{eq:armijo_rule} holds, then the step size is accepted and $h^{(t-1)}=h$, otherwise the step size estimate is reduced to $\alpha h$. This procedure continues until the step size is accepted. The interested reader may refer to  \cite{boyd2004,drummond2018accelerated,nesterov2018lectures,armijo1966minimization} for details.
%
%
\subsection{Second-Order Optimization}
\label{ssec:secondopti}
First-order optimization algorithms such as the one described in Section {\ref{ssec:firstopti}}, while being simple and ensuring convergence to local optima, suffer from a number of disadvantages. One of them is their inability to differentiate between local and global optima; this also affects some second-order optimization methods like Newton's method. Moreover, first-order methods commonly need a large number of iterations to achieve convergence and they also cannot differentiate saddle points---points satisfying first-order optimality conditions but not second-order conditions---from local optima. These disadvantages, however, can be eliminated by using second-order methods. Additionally, even if the
steepest descent (ascent) algorithm is not exactly at a saddle
point (which is unlikely in practice), it
behaves poorly when near them, since gradient steps are small
in directions that escape saddle points (see {\cite{dauphin2014identifying}}).

We can not expect the considered objetive functions to be convex on the design variables; therefore, a modified version of Newton's method for non-convex functions is considered \cite{Paternain2019, dauphin2014identifying}. This method replaces the usual Newton step $\Delta \bm{x}:=H(\bm{x})^{-1}\nabla f(\bm{x})$ with the modified step $\Delta \bm{x}:=\vert H(\bm{x})\vert^{-1}\nabla f(\bm{x})$, where $H(\bm{x})$ is the Hessian matrix of $f(\bm{x})$, and $\vert H(\bm{x})\vert$ is the matrix resulting from taking the spectral decomposition of $H(\bm{x})$ and replacing all its negative eigenvalues by their absolute value. The modified Newton iteration is therefore given as
\begin{align}
	\bm{x}^{(t)}=\bm{x}^{(t-1)}-h^{(t-1)}{\vert H(\bm{x}^{(t-1)})\vert^{-1}}\nabla f(\bm{x}^{(t-1)}).\label{eq:newtonstep}
\end{align}

The modified step not only ensures a faster escape from saddle points than the usual gradient step, but
it also yields quadratic convergence in the direction of eigenvectors associated with positive
eigenvalues. To make this point clearer, let $Q$ be the matrix with orthonormal eigenvectors of the Hessian matrix $H$ as columns
and let $Q_+$ and $Q_-$ be the matrices with eigenvectors of $H$ associated to positive and
negative eigenvalues, respectively. Consider then the projection of the gradient on the spaces spanned by the columns of $Q_+$ and $Q_-$ (hereafter, positive and negative eigenvectors)
$\nabla f_+$ and $\nabla f_-$, respectively. 
The modified Newton step differs from the usual Newton
step only in the direction it takes along the space spanned by the negative eigenvectors, since only negative eigenvalues have their sign changed. As such, quadratic convergence along the positive direction $\nabla f_+$ is expected, while the modification in the negative direction $\nabla f_-$ allows the method to escape saddle
points at an accelerated rate
(see Eq.~(8) and Theorem 3.2 in {\cite{Paternain2019}}, as well as the accompanying discussion).

A variant of the algorithm recalculates $\vert H\vert$ every $m\in\IN$ iterations, i.e.,
\begin{align}
	\bm{x}^{(t\cdot m+i)}=\bm{x}^{(t\cdot m + i - 1)}-h^{(t\cdot m +i-1)}\vert H(\bm{x}^{(t\cdot m)})\vert^{-1}\nabla f(\bm{x}^{(t\cdot m+i-1)}),\label{eq:newtonstepm}
\end{align}
for $i=\{1,\hdots,m\}$ and $k\in\IN_0$, which proves to be advantageous to reduce the computational burden associated with computing the Hessian matrix and its inverse (see Algorithm 1 in {\cite{dauphin2014identifying}}). In our analysis below, both alternatives are implemented and compared. Moreover, the step size for both alternatives
is found through backtracking and the Armijo-Goldstein conditions in
\eqref{eq:armijo_rule} with $\bm{p}^{(t-1)}=\vert H(\bm{x}^{(t-1)})\vert^{-1}\nabla f(\bm{x}^{(t-1)})$.

As for first-order derivatives, the adjoint method allows for the approximation of the far-field second derivatives through solving $2N$ additional integral equations. Accordingly, second-order derivatives for the diffraction efficiency are computed as
\begin{align*}
	\begin{aligned}
		&e_n''(\Gamma;a_1;a_2)=\frac{2k_{y,n}}{k_y}\left(\text{Re}(\cF_n'(\Gamma;a_1))\text{Re}(\cF_n'(\Gamma;a_2))\right.\\
		&+\text{Im}(\cF_n'(\Gamma;a_1))\text{Im}(\cF_n'(\Gamma;a_2))+\text{Re}(\cF_n(\Gamma))\text{Re}(\cF_n''(\Gamma;a_1;a_2))\\
		&+\left.\text{Im}(\cF_n(\Gamma))\text{Im}(\cF_n''(\Gamma;a_1;a_2))\right).
	\end{aligned}
\end{align*}
Though not reported here in detail, both approaches for computing the first and second-order derivatives of the far-field and the efficiency were validated by comparing them
to a finite difference approach.
%
%
\subsection{A Quasi-Newton Method: The BFGS Algorithm}
Finally, we consider another variation on Newton's method in the form of a quasi-Newton update (see Chapter 2.6 in \cite{Gstrang:2007} for other variations) and compare its performance to both steepest descent and modified Newton methods. The BFGS algorithm---originally developed in \cite{broyden1970convergence,fletcher1970new,goldfarb1970family,shanno1970conditioning}---employs the following iteration for the minimization of $f(\bm{x})$: 
\begin{align}
\bm{x}^{(t)}= \bm{x}^{(t-1)}-h^{(t-1)}B^{(t-1)}\nabla f(\bm{x}^{(t-1)}),\label{eq:BFGS}
\end{align}
where the matrix $B^{(t)}$ is obtained from $B^{(t-1)}$ as
\begin{align*}
B^{(t)}=B^{(t-1)}-\frac{B^{(t-1)}s^{(t-1)}{s^{(t-1)}}^{\top}B^{(t-1)}}{{s^{(t-1)}}^{\top}B^{(t-1)}s^{(t-1)}}-\frac{y^{(t-1)}{y^{(t-1)}}^{\top}}{{s^{(t-1)}}^{\top}{y^{(t-1)}}}
\end{align*}
and
\begin{align*}
 s^{(t-1)}=h^{(t-1)}B^{(t-1)}\nabla f(\bm{x}^{(t-1)}),\\
y^{(t-1)}=\nabla f(\bm{x}^{(t)})-\nabla f(\bm{x}^{(t-1)}).
\end{align*}
The matrix $B^{(t)}$ must be positive definite and symmetric for all iterations. This condition is guaranteed if the matrix $B^{(0)}$ is initialized as such and each step size $h^{(t)}$ is chosen so that Wolfe conditions are satisfied, i.e., \eqref{eq:armijo_rule} with $\bm{p}^{(t)}=B^{(t)}\nabla f(\bm{x}^{(t)})$ and
\begin{align}
\bm{p}^{(t-1)}\cdot\nabla f(\bm{x}^{(t-1)}-h^{(t-1)}\bm{p}^{(t-1)})
\leq \gamma \bm{p}^{(t-1)}\cdot\nabla f(\bm{x}^{(t-1)}),
\label{eq:wolfe_curv}
\end{align}
for $\beta\leq\gamma<1$. Though convergence of the algorithm for convex functions is guaranteed, the algorithm, with inexact (Wolfe) line search of the step size, may not converge for non-convex functions 
(see {\cite{dai2002convergence}} and references within).
%
%
\section{Numerical Examples}
\label{sec:numexam}
We present several examples comparing the proposed optimization methods. A standard Galerkin formulation with piecewise polynomials of degree one (previously used in \cite{SilvaOelker2018}) was employed to numerically solve the required BIEs. All computations were performed on a
AMD Opteron 6386 SE server, where parallelization was only used to assemble the relevant matrices. 

The given examples deal with the following optimization problems:
$$\min\limits_\Gamma (e_n(\Gamma)-e_n^{\text{obj}})\quad \text{and} \quad \max\limits_\Gamma e_n(\Gamma),$$
i.e., we seek for a grating geometry that either attains a certain objective diffraction efficiency or maximizes the diffraction efficiency on a given diffraction mode given wave numbers and incidence angles.

The grating period and its relation to the wavelength also plays an important role. Therefore, it may be chosen according to a desired application, for example, considering the Littrow configuration (as considered in Section \ref{ssec:app}) or second Bragg incidence. After presenting these examples, we compare the performances of the different considered algorithms
and, finally, focus on a practical example.
%
%
\begin{figure*}[t!]
	\begin{center}
        \includegraphics[scale=0.5]{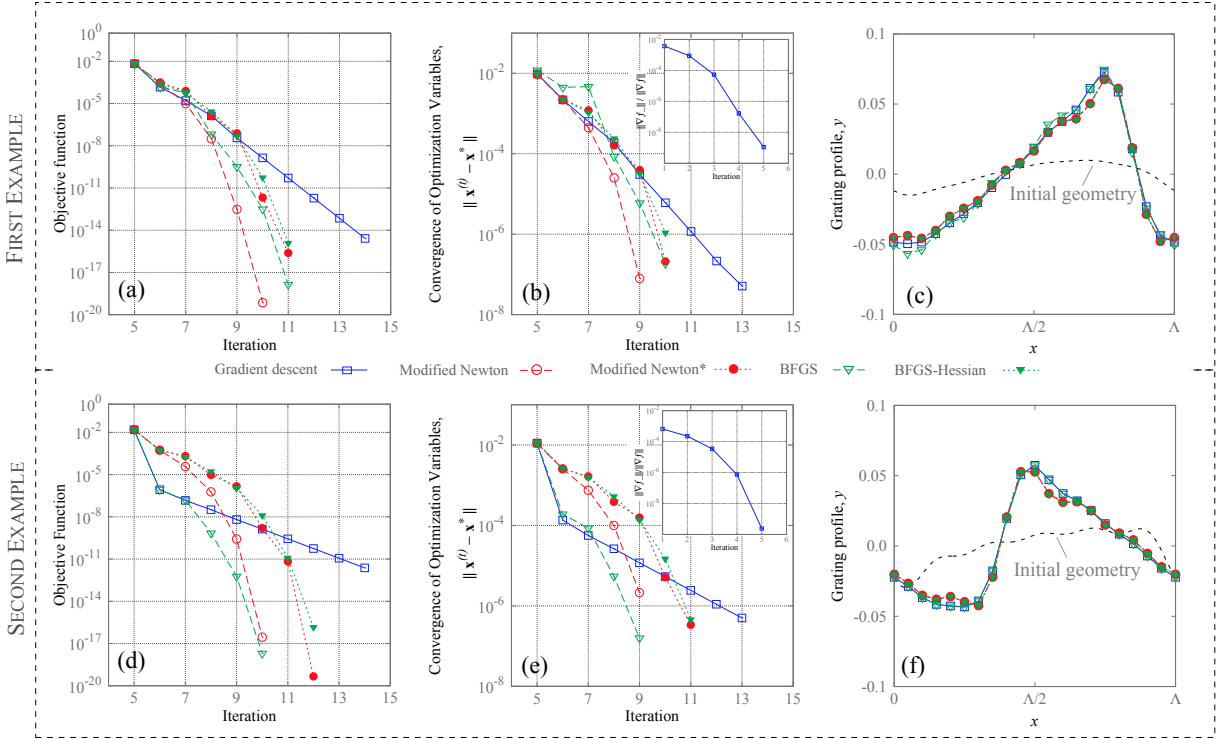}
        \caption{Objective function, error of optimization variables, and optimal geometries are shown
        considering two minimization examples. Top row corresponds to the
		grating profile optimization for the diffraction efficiency in the
		$n=-1$ diffraction order to a target efficiency of $60 \%$ for an
		incidence angle of $\theta_{i} = \pi/6$, and wavenumber $k=40$;
		the bottom row of figures corresponds to the optimization of the diffraction efficiency in the $n=1$ diffraction
		order to an objective efficiency of $65\%$ for an incidence angle of
		$\theta_{i} = \pi/4$, and wavenumber $k=30$.\newline
		The algorithms being compared are specified in the label,
		where the gradient descent algorithm corresponds to \eqref{eq:graddescent}, the
		"Modified Newton" and "Modified Newton$^*$" labels correspond to
		the algorithms in \eqref{eq:newtonstep} and \eqref{eq:newtonstepm}
		with $m=2$, respectively, and the "BFGS" and "BFGS-Hessian" labels correspond
		to the BFGS algorithm with different initializations for $B^{(0)}$ in \eqref{eq:BFGS} -- $B^{(0)}=\mathsf{Id}$ (the identity matrix) and  $B^{(0)}=\vert H\vert$, respectively. Only the steps taken before the initial
		five iterations of gradient descent are shown. For this example, 10 optimization variables were considered ($N=5$).}
				\label{fig:convergence}

	\end{center}        
\end{figure*}
%
%
\subsection{Attaining a Specified Diffraction Efficiency}
\label{ssec:conv}
We start by minimizing an objective function of the form
\begin{align}\label{eq:convfobj}
	\tilde{f}(e_n(\Gamma))=\left(e_n(\Gamma)-e_n^\text{obj}\right)^2,
\end{align}
for different $n\in\IN$ and $e_n^{\text{obj}}\in [0,1]$. Convergence of this objective function and design variables are displayed in Figs. \ref{fig:convergence}. Here, a step size $h=1$ was suitable for our modified Newton method with no backtracking of $h$ required. We take a number of steps of gradient descent from a randomized initial geometry (with a fixed small step) and present  examples for which the algorithms find optimal geometries, $\Gamma_{\text{obj}}$, such that 
$$\tilde{f}(e_n(\Gamma_\text{obj}))=0.$$
Examples when this fails to happen behave similarly as those presented in the following subsection. Specifically, if the initial efficiency is lower than the objective efficiency and the algorithm fails to find a geometry $\Gamma$ such that
$e_n(\Gamma)=e_n^{\text{obj}}$, the minimization of $(e_n(\Gamma)-e_n^{\text{obj}})^2$ behaves like the maximization of $e_n(\Gamma)$.

The examples here consider optimization for a previously defined efficiency as required in applications such as hyperspectral imaging, where some practical parameters yield a required diffraction efficiency {\cite{sabushimike:2018}}. Our first example optimizes the $n=-1$ diffraction order to an objective efficiency of $60\%$ for an incident angle of $\theta_i=\pi/6$ and wavenumber $k=30$. The second example optimizes the $n=1$ diffraction order to an objective efficiency of $65\%$ for an incident angle of $\theta_i=\pi/4$ and wavenumber $k=40$.
The specified incident angles (different than normal incidence) and wavelengths (smaller than the grating period, $\lambda/\Lambda \leq 0.2$) were chosen to show the robustness of the method for scattering problems. For the step size backtracking, the chosen parameters were $\alpha=0.5$ and $\beta=0.2$ for all algorithms to ensure a fair comparison. For the Wolfe line search in \eqref{eq:wolfe_curv}, we chose $\gamma=0.8$.

As indicated in Fig. \ref{fig:convergence}, all implementations of the modified Newton method and the BFGS algorithm---with different initializations for $B^{(0)}$---appear to be superior in rate to the standard gradient descent algorithm. The convergence rate, $q$, of the series of steps $\{\bm{x}^{t}\}_{t\in\IN}$ may be estimated through the formula
\begin{align*}
	q\approx{\log\left(\frac{\norm{\bm{x}^{(t+3)}-\bm{x}^{(t+2)}}}{\norm{\bm{x}^{(t+2)}-\bm{x}^{(t+1)}}}\right)}\left({\log\left(\frac{\norm{\bm{x}^{(t+2)}-x^{(t+1)}}}{\norm{\bm{x}^{(t+1)}-\bm{x}^{(t)}}}\right)}\right)^{-1}.
\end{align*}

This formula
yields an expected convergence rate of $q=1$ for gradient descent,
quadratic convergence for the modified Newton algorithm ($q\approx 2$) and 
super-linear convergence for the BFGS algorithm ($1<q<2$).

Furthermore, in both examples, the gradient descent algorithm required to solve \eqref{eq:BIE} a larger number of times per iteration than either the modified Newton algorithm or BFGS (on both implementations of each algorithm) to find a step size satisfying the corresponding conditions. Indeed, the modified Newton algorithm and the BFGS algorithm with $B^{(0)}=\vert H\vert$ required no backtracking of the step size, while the BFGS algorithm with $B^{(0)}=\mathsf{Id}$ (the identity matrix) required backtracking of the step size.

Quadratic convergence for the modified Newton algorithm and super-linear convergence of the BFGS algorithm are only ensured for convex objective
functions. To explain the observed rates of convergence,
we examine $\nabla f_+$ and $\nabla f_-$. The insets in Figs.~{\ref{fig:convergence}}(b)--(e) display the convergence to zero of $\Vert\nabla f_-\Vert/\Vert\nabla f\Vert$, representing the
fraction of the gradient $\nabla f$ pointing in directions that oppose 
quadratic descent. Thus, the contribution of $\nabla f_-$ is negligible
compared to that of $\nabla f_+$, which accounts for the seemingly
convex behavior of the objective function near its minimum (see {\cite{Paternain2019}}).

\subsection{Efficiency Maximization}
\label{ssec:maxim}
We now focus on maximizing objective functions of the form
\begin{align}
	\tilde{f}(e_n(\Gamma))=e_n(\Gamma)
\end{align}
for different $n\in\IN$ and compare, as before, two different
 examples.
 
 From
 \cite{Paternain2019,dauphin2014identifying}, we expect the
  iteration count of the modified Newton method to be lower than
   that for the first-order method,
as observed in the previous example.
As before, we compare two versions of the modified Newton method (those in \eqref{eq:newtonstep} and \eqref{eq:newtonstepm} with $m=2$) to the first-order method and two instances of the BFGS algorithm with
different initializations for $B^{(0)}$. The algorithms stop if either the gradient of the objective function or the step size at any given iteration fall below a certain tolerance.

\begin{figure}[ht!]
	\centering\includegraphics[width=7cm]{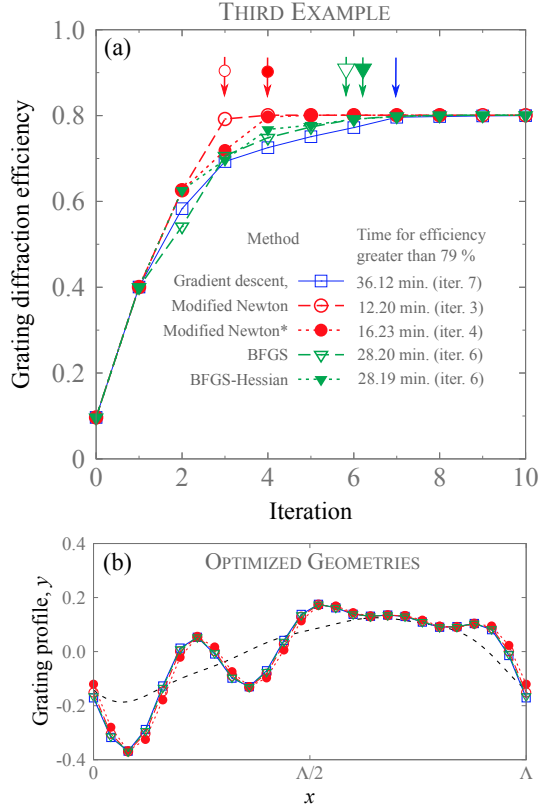}
\caption{{(a) Grating diffraction maximization and (b) optimized geometries are shown for all the studied optimization methods. The parameters correspond to an incidence angle of $\theta = \frac{5}{36}\pi$ and wavenumber $k = 20$ for the maximization of the diffraction efficiency in the $n = 1$ diffraction order, and $N = 4$ (i.e., $8$ optimization variables).}}
	\label{fig:maximization_2}
\end{figure}

\begin{figure}[ht!]
	\centering\includegraphics[width=7cm]{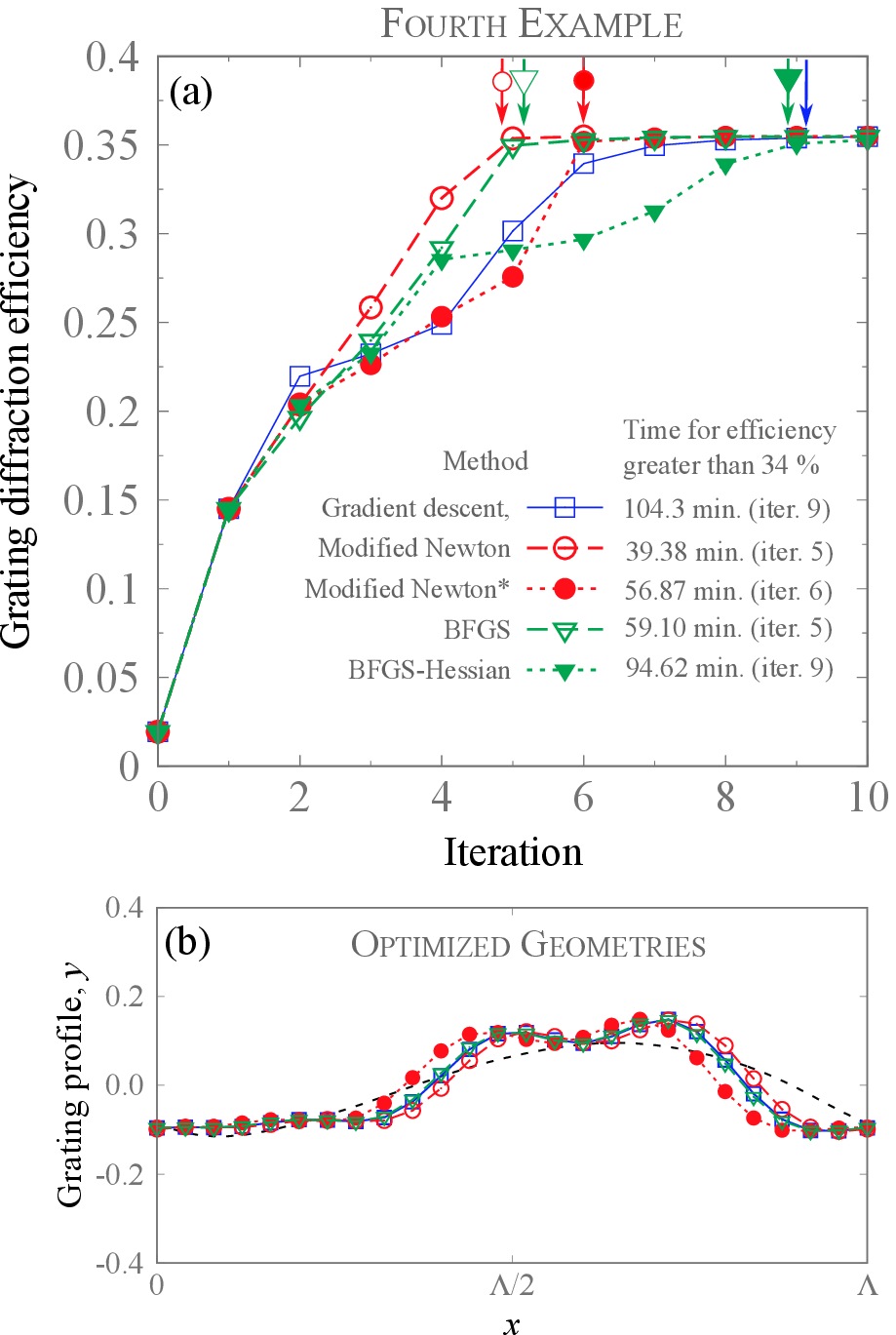}
	\caption{{(a) Grating diffraction maximization and (b) optimized geometries are shown for all the studied optimization methods. The parameters correspond to an incidence angle of $\theta = \frac{1}{4}\pi$ and wavenumber $k = 30$ for the maximization of the diffraction efficiency in the $n = -1$ diffraction order, and $N = 4$ (i.e., $8$ optimization variables).} }
	\label{fig:maximization_3}
\end{figure}

Figures {\ref{fig:maximization_2}} and {\ref{fig:maximization_3}} display the diffraction efficiency, at each iteration, of the target mode being optimized for all three methods, the time each method took to arrive at a tolerance $\varepsilon_{tol}=10^{-2}$ from the maximum efficiency attained, and the final optimized geometries. In general, more than one computation of the objective function was required for the step size backtracking, so the time each method takes to achieve its optimum is not strictly proportional to the number of iterations.

Both versions of the modified Newton method maximize the objective function in fewer 
iterations and time than gradient ascent, as expected, while also converging more rapidly than the BFGS algorithm. 
We note that, in this case, neither quadratic nor super-linear convergence was observed.

%
%

\subsection{Application: Grating Design}
\label{ssec:app}
In this section, gratings are designed and optimized considering the modified Newton algorithm described above. Two grating profiles given by the
linear combination of five sines and cosines ($N = 5$) with different initial conditions are considered and the diffraction efficiency on the first-order mode ($n=1$) is maximized for a wavelength $\lambda_{0} = 300$ nm.

For high diffraction efficiency, echelle gratings are commonly used. For this reason, designed gratings are compared to an echelle grating with a blazed angle of 5.2$^{\circ}$ in a Littrow configuration ($\theta_{inc} = \theta_{n=1}$). Accordingly, through the grating equation, a period of $\Lambda = 1667$ nm is considered 
for the designed gratings, as well as for the echelle grating.

Though a server was used to obtain all previously presented results
at an enhanced level of precision, the studied algorithms may be executed
on regular laptops.
To demonstrate this, the proposed designs
displayed on this section were obtained by running the modified Newton
algorithm on a laptop equipped with an Intel(R) Core(TM) i5-3210M CPU
and 8 GB DDR3 ram memory, taking no longer than $10$ minutes
to converge.

\begin{figure}[ht!]
	\centering\includegraphics[width=7cm]{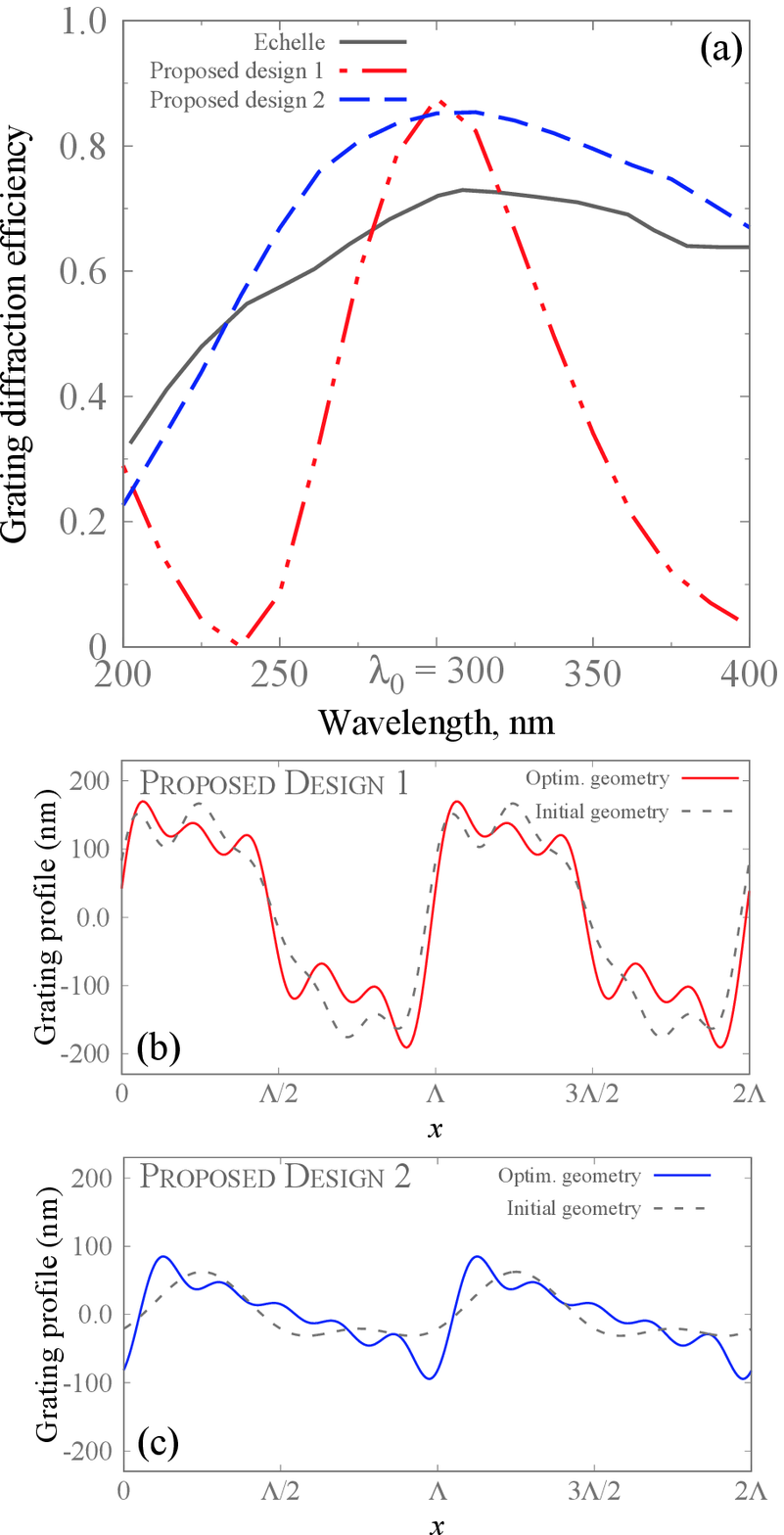}
	\caption{(a) Grating diffraction efficiency for the echelle grating (Newport 53-101R) and two optimized designs, considering $10$ optimization variables ($N=5$). (b) Optimized and initial  geometries for the first design and (c) optimized and initial geometries for the second design.\newline Geometries are optimized for maximal diffraction efficiency on the first-order diffraction mode ($n=1$)  at $\lambda_{0} = 300$ nm in the Littrow configuration.}
	\label{fig:application}
\end{figure}

Figure {\ref{fig:application}} shows the efficiency around $\lambda_{0}$ for both designs and the echelle grating (Newport 53-101R {\cite{Palmer2014}}). The first design gives an efficiency higher than the echelle (72.9 \%) only close to $\lambda_{0}$, with a marked peak of 87.6 \%.
The second design exhibits higher efficiency almost for the entire plotted wavelength range, with an efficiency of 85.2 \% at $\lambda_{0}$.
The sensitivity of the designs to variation was also assessed.
Varying the grating parameters (independently) by 5 \% was
found to reduce the efficiency at $\lambda_{0}$ to 84.8 \%
and 83.4 \% (worst case) for the first and second
designs, respectively.

%
%
\section{Conclusions}
\label{sec:conc}
This work analyzed and compared various optimization algorithms to derive optimal shapes of PEC periodic gratings for high diffraction efficiency for TE polarization. 

The modified Newton method for non-convex functions proved to: (i) achieve quadratic convergence of the optimization variables for the objective function in \eqref{eq:convfobj}; and, (ii) converge to optimal geometries in a lower iteration count than first-order methods for the maximization of grating efficiencies, even if second-order convergence rate is not observed.
The BFGS algorithm outperformed the gradient descent method
for all considered examples, achieving super-linear convergence for the
objective function in \eqref{eq:convfobj}.

Furthermore, the observed convergence rates in the first two numerical examples seem to be
a consequence of the convex nature of the objective function in \eqref{eq:convfobj} (suggested by the low contribution of $\nabla f_-$ to $\nabla f$ displayed in the insets of Figs. {\ref{fig:convergence}} (b) and (e)), hinting that---when such behavior is expected---the modified Newton method should prove advantageous. In particular, applications of inverse problems employ objective functions that should behave in a similar manner near their optima (see \cite{bao1998modeling} and references therein). Moreover,
objective functions similar to those in \eqref{eq:convfobj} may prove useful
for designing gratings with user-chosen efficiency curves \cite{sabushimike:2018}.

The variant of the modified Newton method that recomputes $\vert H \vert$ every $m>1$ iterations (rather than at each iteration), while increasing the number of steps required to converge, proved to be more efficient than its counterpart ($m=1$), when high precision is of interest ($\varepsilon \leq 10^{-3}$). Both variants of the modified Newton method outperformed gradient descent (ascent) in all examples, while
also outperforming the BFGS algorithm on the third and fourth examples.

Finally, two grating profiles were designed using the modified Newton
algorithm and compared to a commercial echelle grating. 
Though both designs differ in their initial geometry,
both outperform the commercial grating at the
target wavelength, but differ significantly in their behavior at neighboring 
wavelengths, revealing the role that starting points play.

The results presented throughout demonstrate the applicability of the
studied algorithms, paving future research 
in other settings, such as transmission problems (dielectrics and semiconductors), where different boundary conditions have to be considered.
The extension to gratings with corners is also of interest for their applications.
However, their optimization presents further challenges, mainly because of the computation of shape derivatives in low-regularity profiles.

\section*{Funding}
This work was sup\hyp{}ported, in part, by Conicyt-PFCHA/Doctorado Nacional/2017-21171791. C.~\hyp{}Jerez-Hanckes work was sponsored by Fondecyt Regular 1171491.

\section*{Disclosures}
The authors declare no conflicts of interest.

\bibliography{references.bib}
\bibliographystyle{plain}

\end{document}